\documentclass[12pt,twoside]{article}
\long\def\symbolfootnote[#1]#2{\begingroup%
\def\thefootnote{\fnsymbol{footnote}}\footnote[#1]{#2}\endgroup}

\usepackage{amsfonts, amssymb, amsmath, amsthm}

\usepackage[colorlinks=true, pdfstartview=FitV, linkcolor=blue,
            citecolor=blue, urlcolor=blue]{hyperref}
\usepackage[usenames]{color}
\definecolor{Red}{rgb}{0.7,0,0.1}
\definecolor{Green}{rgb}{0,0.7,0}
\usepackage{accents}
\usepackage{comment}

\newcommand{\dom}{\mathcal M}

\title
{The    Existence   of Strong  Solutions     to the $3D$ Zakharov-Kuznestov Equation in a Bounded Domain}
\author{Chuntian Wang
}

\date{}

\numberwithin{equation}{section}

\newtheorem{thm}{Theorem}[section]
\newtheorem{lem}{Lemma}[section]
\newtheorem{prop}{Proposition}[section]

\newtheorem{rem}{Remark}[section]

\begin{document}


\maketitle

\vskip-4mm

\centerline{\footnotesize{\it  Department of Mathematics and The Institute for Scientific Computing and Applied Mathematics    }}
\vskip-1mm
\centerline{\footnotesize{\it Indiana University, Bloomington, IN 47405}}

\vskip 0mm
\centerline{\footnotesize{\it \,\, email: \url{wang211@umail.iu.edu}}}

%
%
%

\begin{center}
\large
\date{\today}
\end{center}

\vskip4mm

\tableofcontents

\newpage

\begin{abstract}
 We consider  the Zakharov-Kuznestov (ZK) equation posed in a limited domain $\dom=(0,1)_{x}\times(-\pi /2, \pi /2)^d,$ $ d=1,2$ supplemented with suitable boundary conditions. We  prove that there exists a  solution $u \in \mathcal C ([0, T]; H^1(\dom)) $ to the initial and boundary value problem for the ZK equation in both dimensions $2$ and $3$ for every $T>0$.
 To the best of our knowledge, this is the first result of the global  existence of strong solutions for the ZK equation in   $3D$.

 More importantly, the idea behind the application of    anisotropic estimation to cancel the nonlinear term, we believe, is not only suited for this model but can also be applied to other   nonlinear  equations with similar structures.

 At the same time,  the uniqueness  of  solutions is still open in $2D$ and $3D$ due to the  partially hyperbolic feature of the model.

\end{abstract}
{\noindent \small
 {\it \bf Keywords: Zakharov-Kuznetsov equation,    Korteweg-de Vries equation\\}


\section{Introduction}
\label{sec:introduction}
The Zakharov-Kuznestov (ZK) equation
\begin{equation}\label{eq0r}
\frac{\partial u}{\partial t} + \Delta \frac{\partial  u}{ \partial x  } +  c \frac{  \partial u  }{ \partial x  }+ u \frac{  \partial u  }{ \partial x  }   =f,
\end{equation}
where $u=u(x,  x^\perp, t)$,  $x^\perp =y$ or $x^\perp=(y,z)$,
describes the propagation of nonlinear ionic-sonic waves in a plasma submitted to a magnetic field directed along the $x$-axis.  Here $c>0$ is the sound velocity.  It has been derived formally in a long wave, weakly nonlinear regime from  the Euler-Poisson system in \cite{ZK} and
 \cite {LSp}. A rigorous derivation is provided in \cite{LLS}.  For more general physical references, see   \cite{bonaparitchard1} and \cite{bonapritchatd2}.
When $u$ depends only on $x$ and $t$, (\ref{eq0r}) reduces to the classical Korteweg-de Vries (KdV) equation.

Recently the ZK equation has caught much attention,  not only because it is closely related with the physical phenomena but also because it is the start to explore more general problems that are partly hyperbolic (such as the inviscid primitive equations).



%

Concerning the initial and  boundary value problems of  the Korteweg-de Vries equation posed on a bounded interval $(0,L)$, we refer the interested readers to e.g. \cite{bonasunzhang}, \cite{colinghida}, \cite{tz} and \cite{colingisclon}.

The initial and boundary value problem associated  with (\ref{eq0r}) has been studied
 in the half space $\{ (x, y):\, x>0 \}$ (\cite{Fam2}), on a strip like $\{(x, y):\, x \in \mathbb R, \, 0<y <L \}$ (\cite{BaykovaFaminskii}) or $\{(x, x^\perp):\, 0<x<1, \,x^\perp \in \mathbb R^d, \, d=1,2 \}$ (\cite{Fam3} and \cite{SautTemam}),
  and in a rectangle $ \{ (x, x^\perp): \, 0<x<1, \, x^\perp \in (-\pi/2, \pi/2)^d, \, d = 1, 2 \}$  (\cite{SautTemamChuntian}). Specifically in \cite{SautTemamChuntian},  the authors  have established, for  arbitrary large initial data, the existence of global weak solutions in space dimensions  $2$ and $3$ ($d=1$ and $ 2$ respectively) and a result of uniqueness of such solutions in the two-dimensional case.

As for the existence of  strong solutions,   the global existence  in space dimension $2$ has been proven in a half strip $\{ (x, y):\, x>0,\, y \in (0, L) \}$ in \cite{LT}.  The existence and  exponential decay of regular solutions to the linearized ZK equation in a rectangle $ \{ (x, y): \, x \in (0, L), \, y \in (0, B) \}$ has been studied in \cite{DoroninLarkin}. The local existence of strong solutions in space dimensions $2$ and $3$ is established in \cite{CW}.
In these previous works, the boundary conditions on $x=0,1$   are assumed to be $u\big|_{x=0} = u\big|_{x=1}=u_x \big|_{x=1} =0$; however here we suppose different boundary conditions to serve our purposes.


To the best of our knowledge,   the global existence and uniqueness of regular solutions in $3D$ is still an open problem. In this article, we prove that there exists a global solution $u \in \mathcal C ([0, T]; L^2 (\dom))$ for the initial and boundary value problem of the ZK equation in both $2D$ and $3D$, which we believe, will lead to the global well-posedness of strong of solutions in $3D$ eventually.
 It is interesting to observe that, for the $3D$ ZK equation, the nonlinear term has the same structure as the nonlinear term in the $3D$ Navier-Stkoes equations and that the basic a priori estimates ($L^\infty (0, T; L^2 (\dom))$ and $L^2 (0, T; H^1 (\dom))$) are the same, although the structure of the linear operator is totally different (e.g.    not   coercive  as in (\ref{au}) below).

For the proof we   use the parabolic regularization as in  \cite{SautTemam}, \cite{SautTemamChuntian} and \cite{CW}.
There are four main difficulties. Firstly, as in the case of $3D$   Navier-Stokes equation,
%
the nonlinear term will pose a problem when we apply the Sobolev imbedding in $3D$.
 Secondly,
since the linear operator is not coercive, the $L^p$ estimations (see e.g. \cite{CaoTiti}) does not work. Thirdly,   some assumption on the trace   $u_{xx}\big|^{x=1}_{x=0}$   is necessary for  the estimate  of  $\nabla u\in L^\infty (0, T; L^2 (\dom))$. Finally, to pass to the limit on the boundary conditions, the methods in \cite{SautTemam} and \cite{SautTemamChuntian} are not applicable any more because of the change of the boundary conditions.

To overcome these   difficulties, firstly we utilize the anisotropic resonance of the  term $u_{xxx}$ and the nonlinear term $uu_x$   to cancel  $uu_x$, which leads to a  bound of the $H^1$ norm over $(0, T)$ for $u$.
This step  of canceling the nonlinear term may also be applied  to other nonlinear   equations with similar structures.
Next, we suppose periodic boundary conditions of $u$ and $u _{x^j}$ at $x=0$, $1$, $j=1, 2$, so that  the trace $u_{xx}\big|^{x=1}_{x=0}$ now vanishes. Finally, we investigate a bound independent of $\epsilon$ for $u^\epsilon_{xxx}$ in $L^{3/2} (I_x; Y)$, with $Y$  a Banach space in $x^\perp$ and $t$, which facilitates the passage to the limit on the traces of  $u _{x^j}$ at $x=0$, $1$, $j=1, 2$.

 However    the uniqueness of solutions is still open in both $2D$ and $3D$, even with such a regularity and all the periodic boundary conditions satisfied. In particular,  the methods   in \cite{SautTemam} and \cite{SautTemamChuntian} can not be adapted to our case due to the lack of the boundary condition $u_x=0$ at $x=1$.



The  article is organized as follows. Firstly we introduce the basic settings of the equation in Section \ref{sec0}.
Secondly we introduce the  parabolic regularization  as in \cite {SautTemam} and \cite{SautTemamChuntian} (Section \ref{regu}). Then we derive the
  estimates independent of $\epsilon$ for $u^\epsilon$ in $L^\infty(0, T; L^2 (\dom))$ (Section \ref{l2e}),  $\nabla u^\epsilon$ in $L^\infty(0, T; L^2 (\dom))$ (Section \ref{h1e})  and for  $  u^\epsilon _{xxx}$ in $L^{3/2}(I_x;  H^{-1}_t (0, T;\, H^{-4}( I_{x^\perp}   )) )$    (Section \ref{other}).
Eventually we can pass to the limit on the parabolic regularization and  the traces and deduce the global existence of solutions $u\in \mathcal C ([0, T]; \, H^1 (\dom))$   (Section \ref{gob}). Finally,    we    discuss about   the difficulties in the attempt of proving the uniqueness of solutions (Section \ref{unique}).


\section{ZK equation in a rectangle in dimensions $2$ and $3$}
\label{sec0}
We aim to  study the ZK equation:
\begin{equation}\label{eq01}
\frac{\partial u}{\partial t} + \Delta \frac{\partial  u}{ \partial x  } +  c \frac{  \partial u  }{ \partial x  }+ u \frac{  \partial u  }{ \partial x  }   =f,
\end{equation}
in a rectangle or  parallelepiped domain  in $\mathbb R^n$  with $n= 2$ or $ 3$, denoted as  $\dom=(0,1)_{x}\times(-\pi /2, \pi /2)^d$,   with $ d=1 $ or $2$, $\Delta u=u_{xx}+\Delta^\perp u$, $\Delta^\perp u=u_{yy}$ or $u_{yy}+u_{zz}$ depending on the dimension.
In the sequel we will use the notations $I_{x}=(0,1)_{x},$ $ I_{y}=(- \pi /2,  \pi /2)_{y}$,  $I_{z}=(- \pi /2,  \pi /2)_{z}$, and $I_{x^\perp} =I_y$ or $ I_{y} \times I_{z}$. We assume the boundary conditions of $u$, $u_x$ and $u_{xx}$ on $x=0, 1$ to be periodic:
\begin{equation}\label{dirch}
 u(0, \, x^\perp,\, t)=u(1, \, x^\perp,\, t),
\end{equation}
\begin{equation}\label{eq4}
\begin{split}
 u _{x}(0, \, x^\perp,\, t)=u_{x}(1, \, x^\perp,\, t), \,\,\,\,\,\,\, u _{xx}(0, \, x^\perp,\, t)=u_{xx}(1, \, x^\perp,\, t).
  \end{split}
 \end{equation}


For the boundary conditions in the $y$ and $z$ directions,  we will choose  either the Dirichlet boundary conditions
\begin{equation}\label{eq107}
 u=0 \mbox{ at }  y=\pm\frac{\pi}{2} \mbox { and }z=\pm\frac{\pi}{2} ,
\end{equation}
or the periodic boundary conditions
\begin{equation}\label{eq1-26p}
\begin{split}
 u\big|^ {y= \frac {\pi}{2}}_{y=-\frac {\pi}{2}}=u_y\big|^ {y= \frac {\pi}{2}}_{y=-\frac {\pi}{2}} =0,\\
 u\big|^{z= \frac {\pi}{2}} _{z=-\frac {\pi}{2}}=u_z\big|^{z= \frac {\pi}{2}} _{z=-\frac {\pi}{2}} = 0.
 \end{split}
\end{equation}

The initial condition reads:
\begin{equation}\label{eq33}
 u(x, \, x^\perp,\, 0)=u_0(x, \, x^\perp).
\end{equation}

 We   study the initial and boundary value problem (\ref{eq01})-(\ref{eq4})
  and  (\ref{eq33}) supplemented with the boundary condition (\ref{eq107}), that is,   the  Dirichlet case   on the $x^\perp$ boundaries,  and we will make some remarks on the extension to the periodic boundary condition case.

We  denote by $|\cdot|$ and $(\cdot, \cdot) $ the norm and the inner product of $L^2 (\dom)$, and by $[\cdot]_2$ the following seminorm which will be useful in the sequel:
\begin{equation}\label{seminorm}
\left (\int_{\dom}      u_{xx}  ^2 +  u_{yy}  ^2 +
 u_{yy}  ^2   \,d\, \dom \right ) ^{1/2}=:  [u]_2, \,\,\,\,\,\, u \in H^2 (\dom).
\end{equation}

\section{Existence of solutions $u\in  \mathcal C  ([0, T]; H^1 (\dom))$  in  dimensions $2$ and $3$}
\label{thee}
To prove this result, we   use the parabolic regularization as  in \cite{SautTemamChuntian}, but with   different  boundary conditions.
For the sake of simplicity we   only treat the more complicated case when $d=2$.

\subsection{Parabolic regularization}
\label{regu}

To begin with, we recall the  parabolic regularization introduced in \cite{SautTemam} and  \cite{SautTemamChuntian},
that is, for $\epsilon >0$ ``small'', we consider the parabolic equation,
\begin{eqnarray}\label{10-11}\displaystyle
\begin{cases}
&\dfrac{\partial u^\epsilon}{\partial t}+ \Delta \dfrac{\partial u^\epsilon }{\partial x}+ c \dfrac{\partial u^\epsilon }{\partial x} +     u^\epsilon \dfrac{\partial u^\epsilon }{\partial x}  +  \epsilon L u ^\epsilon    =f,\\
&u^\epsilon (0)=u_0,
\end{cases}
\end{eqnarray}
where
  \begin{equation*}
L u^\epsilon : =  \dfrac{\partial^4 u^\epsilon}{\partial x^4}   +
 \dfrac{\partial^4 u^\epsilon}{\partial y^4}   + \dfrac{\partial^4 u^\epsilon}{\partial z^4} ,
\end{equation*}
  supplemented with the boundary conditions (\ref{dirch})-(\ref{eq107}) and the additional boundary conditions
\begin{equation}\label{uxxx1}
 u^\epsilon _{xxx}(0, \, x^\perp,\, t)=u^\epsilon_{xxx}(1, \, x^\perp,\, t),
\end{equation}
\begin{equation}\label{10-4}
  u^\epsilon_{yy}=0 \mbox{ at } y=\pm\frac{\pi}{2}, \,\,\,   u^\epsilon_{zz}=0 \mbox{ at } z=\pm\frac{\pi}{2}.
\end{equation}
Note that from (\ref{eq4}) and (\ref{uxxx1}) we infer
\begin{equation}\label{uxxx}
 u^\epsilon _{x^j}(0, \, x^\perp,\, t)=u^\epsilon_{x^j}(1, \, x^\perp,\, t),\,\,\,\,j=1 ,\,2,\,\,3.
\end{equation}
We also note that since $u^\epsilon_{yy}\big|^{x=1}_{x=0}= u^\epsilon_{zz}\big|^{x=1}_{x=0} =0$,  (\ref{uxxx}) is equivalent to
\begin{equation}\label{delta}
\Delta u^\epsilon \big|^{x=1}_{x=0}  = 0.
\end{equation}

It is a classical result (see e.g. \cite{JLLions},     \cite{MR0241822} or also  \cite{SautTemamChuntian}) that      there exists a unique   solution to the parabolic problem which is sufficiently regular for all the subsequent calculations to be valid; in particular,
we have
\begin{equation}\label{uep}
u^\epsilon \in   L^2(0, T; H^4(\dom) ) \cap \mathcal C ^1 ([0, T]; H^2 (\dom)) .
\end{equation}

\subsection{Estimates independent of $\epsilon$ }
\label{seca}

We   establish the     estimates independent of $\epsilon$  for various norms of the  solutions.
\subsubsection{$L^2 $ estimate independent of $\epsilon$}
\label{l2e}
We first show a bound independent of $\epsilon$ for $u^\epsilon$ in $L^\infty (0, T; L^2 (\dom))$.

\begin{lem}\label{main0}
We assume that
\begin{equation}\label{initial2}
u_{0} \in L^2 (\dom) ,
\end{equation}
\begin{equation}\label{f}
f \in L^2(0, T; L^2 (\dom)).
\end{equation}
Then  for every $T>0$ the following estimates independent of $\epsilon$ hold:
\begin{equation}\label{set1}
 u^\epsilon \,\,\mbox{is\,\, bounded\,\,in}\,\,L^\infty(0, T; \,L^2(\dom)),
 \end{equation}
 \begin{equation}\label{set1p}
 \sqrt \epsilon \,  u^\epsilon \,\,\mbox{is\,\, bounded\,\,    in}\,\,L^2(0, T; \,H^2  (\dom)).
 \end{equation}

\end{lem}

\noindent\textbf{Proof.}
 As in \cite{SautTemamChuntian}, we multiply (\ref{10-11}) with $u$, integrate over $\dom$ and integrate by parts,     dropping the superscript $\epsilon$ for the moment we find:
\begin{align*}
  & \bullet\,\,  \int_{\dom} \frac{\partial u•}{\partial t} u \,d\,\dom = \frac{1}{2} \frac {d}{dt} |u|^2 ,\quad\quad\quad\quad\quad\quad\quad\quad\quad\quad\quad\quad\quad\quad\quad\quad\quad\quad\quad\quad\quad\quad\quad
  \end{align*}
  \begin{equation}\label{au}
  \begin{split}
\bullet  \int_{\dom}\Delta u_x  &  \, u  \, d\, \dom + \int_{\dom}  cu_x \, u  \, d\, \dom = (\mbox{thanks to (\ref{dirch})}) \quad\quad\quad\quad\quad\quad\quad\quad\quad\quad\\
&=  -   \int_{\dom}  \nabla u _x \, \nabla u \, d\, \dom + \frac{c}{2}    \int_{I_{x^\perp}}   u ^2 \big|^{x=1}_{x=0}\,d\,x^\perp \\
    &  = -\frac{1}{2}  \int_{I_{x^\perp}}  ( \nabla u )^2 \big|^{x=1}_{x=0}    \,dx^\perp   + \frac{c}{2}    \int_{I_{x^\perp}}   u ^2 \big|^{x=1}_{x=0}\,d\,x^\perp\\
    &= (\mbox{thanks to (\ref{dirch})  and  (\ref{uxxx})} ) =0   ,
    \end{split}
    \end{equation}
    \begin{align*}
  & \bullet\,\,  \int_{\dom}  u u_x \, u \, d\, \dom     =    \int_{\dom} \frac{\partial }{\partial x} \left ( \frac{u^3}{3}  \right )\, d\, \dom= (\mbox{thanks to (\ref{dirch})} )   = 0 ,\\
  & \bullet\,\,    \epsilon \int _{\dom} u_{xxxx}\, u \,d\, \dom    = (\mbox{thanks to (\ref{dirch}) and (\ref{uxxx})})\\
  &\quad\quad\quad\quad \quad\quad\quad\quad\quad= - \epsilon \int _{\dom} u_{xxx}\, u_x     \,d\, \dom = (\mbox{thanks to (\ref{uxxx})}) =   \epsilon \int _{\dom}  u_{xx} ^2    \,d\, \dom , \\
     & \bullet\,\,    \epsilon \int _{\dom} \left ( u_{xxxx}+ u_{yyyy} +u_{zzzz} \right)\, u \,d\, \dom  =\epsilon \int_{\dom}      u_{xx}  ^2 +  u_{yy}  ^2 + u_{yy}  ^2   \,d\, \dom \\
& \quad\quad\quad\quad \quad\quad\quad\quad\quad= (\mbox{thanks to }(\ref{seminorm}))= \epsilon [u]^2_2 , \\
  & \bullet\,\,  \int_{\dom} f u \, d\, \dom \leq \frac{1}{2} \left |f \right|^2  + \frac{1}{2}  |u|^2 .\\
\end{align*}
Hence we find
  \begin{equation}\label{eq10-6}
\frac{d}{dt} |u^\epsilon(t)|^2  + 2\epsilon [u^\epsilon]^2_2
  \leq \left |f \right|^2 + |u^\epsilon |^2.
\end{equation}
Using the Gronwall lemma we classically infer
\begin{equation}\label{eq10-9}
\sup_{t\in (0, T)}|u^\epsilon(t)|^2  + \epsilon\, \int^T _0 [u^\epsilon]^2_2 \, dt \leq const:=\mu_1,
\end{equation}
where  $\mu_i$ indicates a constant depending only on the data $u_0$, $f$, etc, whereas $C^\prime $   below is an absolute constant. These constants may be different at each occurrence.
Let us admit for the moment the following:
\begin{lem}\label{h221}
\begin{equation}\label{[]2}
 |u^\epsilon|^2_{H^2 (\dom)}  \leq   C^\prime  \left (  [u^\epsilon]^2_2     + |u^\epsilon|^2\right) .
\end{equation}
\end{lem}
By the previous lemma, we have
\begin{align*}
\epsilon \, \int ^T_0   |u^\epsilon|^2_{H^2 (\dom)} \, dt  &\leq   C^\prime  \left (    \epsilon \, \int ^T_0 [u^\epsilon]^2_2  \, dt  +   \epsilon \, \int ^T_0 |u^\epsilon|^2 \, dt \right)\\
& \leq C^\prime  \left (   \epsilon \, \int ^T_0 [u^\epsilon]^2_2  \, dt  +    \epsilon \,T \,\sup_{t\in (0, T)}|u^\epsilon(t)|^2 \right)\\
&\leq (\mbox{thanks to (\ref{eq10-9})})\\
 & \leq const:=\mu_2,
\end{align*}
which implies  (\ref{set1p}).
Thus Lemma \ref{main0} is proven once we have proven Lemma (\ref{h221}).

\noindent{\textbf{Proof of Lemma \ref{h221}.}} We first observe that  using the generalized Poincar\'e inequality (see  \cite{Temam3})  we have
\begin{equation}\label{poincare}
|u^\epsilon_x- \int^1 _0 u^\epsilon_x \,dx |_{L^2 (I_x)} \leq C^\prime  |u^\epsilon_{xx}| _{L^2 (I_x)}.
\end{equation}
Thanks to  (\ref{dirch}), we have $\int^1 _0 u^\epsilon_x \,dx = u^\epsilon|^{x=1}_{x=0} = 0$, and hence (\ref{poincare}) implies
\begin{equation*}
|u^\epsilon_x   |_{L^2 (I_x)} \leq C^\prime  |u^\epsilon_{xx}| _{L^2 (I_x)}.
\end{equation*}
Squaring both sides and integrating both sides on $I_{x^\perp}$, we find
\begin{equation}\label{poincare1}
|u^\epsilon_x   |  \leq C^\prime  |u^\epsilon_{xx}| .
\end{equation}
Similarly we can show that  $|u^\epsilon_y   |  \leq C^\prime  |u^\epsilon_{yy}| $ and $|u^\epsilon_z   |  \leq C^\prime  |u^\epsilon_{zz}| $, which implies
\begin{equation}\label{poincare2}
|\nabla u^\epsilon   |  \leq C^\prime  [u^\epsilon]_2 .
\end{equation}
Next we see that, for smooth functions
\begin{equation}\label{mix1}
\begin{split}
|u^\epsilon_{xy}|^2 &=  (\mbox {thanks to (\ref{dirch}) and (\ref{uxxx})})\\
& = - \int _{\dom}  \, u^\epsilon_{y}  u^\epsilon_{xxy}    d \dom \\
&  = (\mbox{thanks to (\ref{eq107})})\\
&  =\int _{\dom}  \, u^\epsilon_{yy}  u^\epsilon_{xx}    d \dom\\
& \leq |u^\epsilon_{xx}|^2 + |u^\epsilon_{yy}|^2 \leq [u^\epsilon]_2 ^2.
\end{split}
\end{equation}
Similarly we can prove that $|u^\epsilon_{xz}|  \leq [u^\epsilon]_2$ and $|u^\epsilon_{zy}|  \leq [u^\epsilon]_2$, and hence
\begin{equation}\label{mix}
\begin{split}
|u^\epsilon_{xy}|^2 + |u^\epsilon_{xz}|^2 +    |u^\epsilon_{yz}|^2    \leq C^\prime   [u^\epsilon]_2 ^2.
\end{split}
\end{equation}
Then inequality (\ref{mix1}) and  (\ref{mix})  extend by continuity to all $H^2$ function periodic in $x$ and satisfying  (\ref{eq107})   and (\ref{10-4}).  Finally from   (\ref{mix})   and   (\ref{poincare2}) we deduce (\ref{[]2}). \qed



\subsubsection{$H^1$ estimate independent of $\epsilon$}
\label{h1e}

Now we establish the key observation, a bound independent of $\epsilon$ for  $\nabla u^\epsilon$ in  $L^\infty (0, T; L^2(\dom))$.
\begin{prop}\label{localbounds1}
Under the same assumptions as in Lemma \ref{main0},   we further suppose that
\begin{equation}\label{u02}
u_0 \in H^1 (\dom) \cap L^3(\dom),
\end{equation}
\begin{equation}\label{f1}
f \in L^2(0, T; H^2  (I_x; \,H^2\cap H_0^ 1 (I_{x^\perp}))) \cap L^2 (0, T; L^\infty (\dom)),
 \end{equation}
and  $f$ and $f_x$ assume the periodic boundary conditions on $x=0$, $1$. Then for every $T>0$,   the following   estimates  independent of $\epsilon$ hold:
\begin{equation}\label{main8}
u^\epsilon \,\,\mbox{is\,\, bounded\,\,in}\,\,L^\infty(0, T; \,H^1 (\dom)),
\end{equation}
\begin{equation}\label{main8p}
  \sqrt \epsilon \, \nabla  u^\epsilon_{xx}, \,\,\sqrt \epsilon \, \nabla  u^\epsilon_{yy},\,\,\sqrt \epsilon \, \nabla  u^\epsilon_{zz} \,\,\mbox{are\,\, bounded\,\,    in}\,\,L^2(0, T; L^2 (\dom)).
\end{equation}
\end{prop}
\noindent \textbf{Proof.} We multiply (\ref{10-11}) with $-\Delta u^\epsilon - \frac{1}{2}\left(u^\epsilon \right)^2$, integrate over $\dom$ and integrate by parts.
Firstly we show the calculation details of the multiplication by $\Delta u^\epsilon$, integration over $\dom$ and integration by parts (dropping the super index of $\epsilon$ for the moment):
\begin{align*}
& \bullet \int_{\dom} u_{t} \, \Delta u \, d\dom = -  \int_{\dom}\nabla u_{t} \, \nabla u \, d\dom + \int_{\partial \dom} u_t \, \frac{\partial u}{\partial n}\,d\,  \partial \dom\\
&\,\,\,\,\,\,\,\,\,\quad\quad\quad \,\,\,\,\,\,\,\,\,\,\,\,\,\,= (\mbox{thanks to }(\ref{dirch}) \mbox{ and } (\ref{uxxx}))\\
  &\,\,\,\,\,\,\,\,\,\quad\quad\quad \,\,\,\,\,\,\,\,\,\,\,\,\,\, = -  \int_{\dom}\nabla u_{t} \, \nabla u \, d\dom =  -\frac{1}{2} \frac{d}{dt} \left |  \nabla u \right |^2,\\
&  \bullet \int_{\dom} \Delta u_{x} \, \Delta u\, d\dom = \int_{\dom} \frac{\partial }{\partial x} \left ( \frac{(\Delta u )^2}{2} \right )  \, d\dom = \frac{1}{2} \int_{I_{x^\perp}} (\Delta u )^2 \big|^{x=1} _{x=0} \, d \,  {I_{x^\perp}} = (\mbox{thanks to (\ref{delta})})=0   , \\
&  \bullet c\int_{\dom} u_{x} \, \Delta u\, d\dom =  c\int_{\dom} u_{x} \,  u_{xx} +   u_x\, \Delta^\perp u \, d\dom \\
& \,\,\,\,\,\,\,\,\,\quad\quad\quad \,\,\,\,\,\,\,\,\,\,\,\,\,\,  =(\mbox{thanks to (\ref{dirch})}) =  c\int_{\dom} \frac{\partial}{\partial x} \left (  \frac{\left(u_x\right)^2 }{2}  \right  )  \, d\dom  - c  \int_{\dom} \nabla^\perp u_{x} \, \nabla^\perp  u\, d\dom\\
   & \,\,\,\,\,\,\,\,\,\quad\quad\quad \,\,\,\,\,\,\,\,\,\,\,\,\,\,  = c\int_{\dom} \frac{\partial}{\partial x} \left (  \frac{\left(u_x\right)^2 }{2}  \right  )  \, d\dom    - c  \int_{\dom}  \frac{\partial}{\partial x}\left ( \frac{(\nabla^\perp u)^2}{2}\right )\, d\dom\\
        & \,\,\,\,\,\,\,\,\,\quad\quad\quad \,\,\,\,\,\,\,\,\,\,\,\,\,\, =  (\mbox{thanks to (\ref{uxxx}) and (\ref{dirch})})  =0     ,
        \end{align*}
        \begin{align*}
  & \bullet\,\,      \int _{\dom}   u_{xxxx} \,   u_{xx} \,d\, \dom = (\mbox{thanks to } (\ref{uxxx}) )=  - \int _{\dom}   u^2_{xxx}   \,d\, \dom,\\
  & \bullet\,\,      \int _{\dom}   u_{xxxx} \,   u_{yy} \,d\, \dom =   (\mbox{thanks to    (\ref{dirch})-(\ref{eq107}) and (\ref{uxxx})})    =  - \int _{\dom}   u^2_{xxy}  \,d\, \dom,\\
  & \bullet\,\,      \int _{\dom}   u_{xxxx} \,   u_{zz} \,d\, \dom = (\mbox{thanks to    (\ref{dirch})-(\ref{eq107}) and (\ref{uxxx})})    =  - \int _{\dom}   u^2_{xxz}   \,d\, \dom,\\
  & \bullet\,\,      \int _{\dom}   u_{yyyy} \, \Delta u \,d\, \dom = (\mbox{thanks to } (\ref{eq107}) \mbox{  and  } (\ref{10-4}))=  - \int _{\dom}   u_{yyy} \, \Delta u_y \,d\, \dom\\
  &  \,\,\,\,\,\,\,\,\,\quad\quad\quad \,\,\,\,\,\,\,\,\,\, \quad\quad\quad   =  (\mbox{thanks to } (\ref{10-4}))= \int _{\dom}   u_{yy} \, \Delta u_{yy} \,d\, \dom\\
  &  \,\,\,\,\,\,\,\,\,\quad\quad\quad \,\,\,\,\,\,\,\,\,\,  \quad\quad\quad  =  (\mbox{thanks to } (\ref{10-4}))= - \int _{\dom}  \left (\nabla u_{yy}\right ) ^2  \,d\, \dom,\\
  & \bullet\,\,      \int _{\dom}   u_{zzzz} \, \Delta u \,d\, \dom =  - \int _{\dom}  \left (\nabla u_{zz}\right ) ^2  \,d\, \dom,\\
& \bullet \int_{\dom} f \, \Delta u \, d\dom =(\mbox{thanks to   (\ref{f1})})  =     \int_{\dom}\Delta f  \, u \, d \dom,
\end{align*}
Hence we find after changing the sign,
\begin{equation}\label{key1}
\begin{split}
\frac{1}{2} \frac{d}{dt} \left|  \nabla u^\epsilon \right|^2   -   &\int _{\dom} u^\epsilon u^\epsilon_x \, \Delta u^\epsilon  \, d \dom +\epsilon [\nabla u^\epsilon]^2_2=  -  \int_{\dom}\Delta f  \, u^\epsilon \, d \dom.
\end{split}
\end{equation}

Next we show the calculation details of the multiplication by $\left( u^\epsilon\right)^2 $, integrating over $\dom$ and integrating by parts:
\begin{align*}
& \bullet \int_\dom u_t  u^2 \, d\, \dom = \int_\dom  \frac{\partial }{\partial t} \left (\frac{u^3}{ 3}  \right )  \, d\, \dom    =
\frac{1}{3}  \frac{d}{dt} \left (\int_{\dom} u^3 \, d \dom \right ),\\
&  \bullet  \int_{\dom}   \Delta u_{x} \,  u^2 \, d\dom =- 2 \int_{\dom} \Delta u \, u u_x\, d\dom + \int_{I_{x^\perp}} \Delta u \, u^2 \big|^{x=1}_{x=0} d \,  {I_{x^\perp}}\\
& \,\,\,\,\,\,\,\,\,\,\,\,\,\,\,\quad\quad\quad\quad\quad = (\mbox{thanks to } (\ref{delta}) \mbox{ and } (\ref{dirch}))   = - 2 \int_{\dom} \Delta u \, u u_x\, d\dom,\\
&  \bullet c  \int_{\dom}   u_{x} \,  u^2 \, d\dom =   c  \int_{\dom}  \frac{\partial }{\partial x} \left ( \frac{u^3}{3}  \right) \, d\dom
= \frac{c}{3}\int_{I_{x^\perp}}  u^3 \big|^{x=1}_{x=0}\,  d \,  {I_{x^\perp}} = (\mbox{thanks to } (\ref{dirch}))=0,\\
& \bullet \int_{\dom} uu_x \, u^2 \, \dom =\int_{\dom}   \frac{\partial}{\partial x} \left (\frac { u^4}{4} \right ) d \dom =\frac{1}{4}\int_{I_{x^\perp}}  u^4 \big|^{x=1}_{x=0}\, d \,  {I_{x^\perp}} = (\mbox{thanks to }  (\ref{dirch}))=0,\\
 & \bullet\,\,      \int _{\dom}   u_{xxxx} \, u^2 \,d\, \dom = (\mbox{thanks to } (\ref{dirch})  \mbox{ and }  (\ref{uxxx}))=  - 2 \int _{\dom}   u_{xxx} \,   u_x \,u\,d\, \dom,\\
   & \bullet\,\,      \int _{\dom}   u_{yyyy} \, u^2 \,d\, \dom = (\mbox{thanks to } (\ref{eq107}))=  - 2 \int _{\dom}   u_{yyy} \,   u_y \,u\,d\, \dom,\\
  & \bullet\,\,      \int _{\dom}   u_{zzzz} \, u^2 \,d\, \dom = (\mbox{thanks to } (\ref{eq107}))=  - 2 \int _{\dom}   u_{zzz} \,   u_z \,u\,d\, \dom.
\end{align*}
Hence we find
\begin{equation}\label{key2}
\begin{split}
\frac{1}{3}  \frac{d}{dt}\left ( \int_{\dom} \left (u^\epsilon \right)^3\, d \dom \right )& -2 \int_{\dom} \Delta u^\epsilon \, u^\epsilon u^\epsilon_x\, d\dom =\\
 &2\epsilon \int _{\dom}   u^\epsilon_{xxx} \,   u^\epsilon_x \,u^\epsilon  + u^\epsilon_{yyy} \,   u^\epsilon_y \,u^\epsilon  + u^\epsilon_{zzz} \,   u^\epsilon_z \,u^\epsilon  \, d\, \dom +  \int_{\dom}  f  \left (u^\epsilon\right ) ^2 \, d \dom .
\end{split}
\end{equation}
Adding  (\ref{key1}) to (\ref{key2}) multiplied by $- 1/2$, we observe that the terms $ \int_{\dom} \Delta u^\epsilon \, u^\epsilon u^\epsilon_x\, d\dom$ get canceled,  which yields
\begin{equation*}\label{key3}
\begin{split}
\frac{1}{2} \frac{d}{dt} \left |  \nabla u^\epsilon  \right |^2 + \epsilon [\nabla u^\epsilon]^2_2 = &\frac{1}{6} \, \frac{d}{dt} \left (\int_{\dom} \left (u^\epsilon \right)^3 \, d \dom \right )  \\
 &\,\,\,\,\,\,\,\,\,\, - \epsilon \int _{\dom}   u^\epsilon_{xxx} \,   u^\epsilon_x \,u^\epsilon   + u^\epsilon_{yyy} \,   u^\epsilon_y \,u ^\epsilon + u^\epsilon_{zzz} \,   u^\epsilon_z \,u^\epsilon  \, d\, \dom\\
 & \,\,\,\,\, \,\,\,\, - \int_{\dom} \Delta f   u^\epsilon \, d \dom - \frac{1}{2} \int_{\dom}   f  \left (u^\epsilon\right ) ^2  \, d \dom .
\end{split}
\end{equation*}
Integrating both sides in time from $0$ to $t$, we obtain for every $t\in (0, T)$,
\begin{equation}\label{key4}
\begin{split}
\frac{1}{2} \left |  \nabla u^\epsilon (t) \right|^2 + \epsilon   \int^t _0 [\nabla u^\epsilon]^2_2 \, ds &= \frac{1}{6}    \int_{\dom}  \left (u^\epsilon (t) \right)^3\, d \dom + \kappa_0\\
&\,\,\,\,\, \,\,\,\,-\epsilon\int^t _0 \int _{\dom}   u^\epsilon_{xxx} \,   u^\epsilon_x \,u^\epsilon + u^\epsilon_{yyy} \,   u^\epsilon_y \,u ^\epsilon + u^\epsilon_{zzz} \,   u^\epsilon_z \,u ^\epsilon \, d\dom   \, ds   \\
 & \,\,\,\,\, \,\,\,\,-  \int^t _0  \int_{\dom} \Delta f   u ^\epsilon\, d \dom \, ds    - \frac{1}{2} \int^t _0 \int_{\dom}   f  \left (u^\epsilon\right )^2 \, d \dom\, ds ,
\end{split}
\end{equation}
where
\begin{equation*}\label{kappa0}
\kappa_0 := \frac{1}{2} \left |  \nabla u_0 \right |^2   -\frac{1}{6}     \int_{\dom} u_0^3 \, d \dom  .
\end{equation*}
We estimate each term on the right-hand-side of (\ref{key4}); we will use here the interpolation space $H^{1/2}(\dom)$ as defined in \cite{LionsMagenes1} where it is shown that $ H^{1/2}(\dom)\subset L^3(\dom) $ in dimension $3$ with a continuous embedding.
 Dropping the superscript $\epsilon$   for the moment we then find:
\begin{equation*}\label{u3}
\begin{split}
\left | \frac{1}{6}    \int_{\dom} u^3(t) \, d \dom \right | &\leq   \frac{1}{6}   |u(t)|^3 _{L^3(\dom)} \\
&\leq   C^ \prime  |u(t)|^3 _{H^{1/2}(\dom)} \\
 & \leq C^\prime |u(t)|^{3/2} \left |\nabla u (t) \right| ^ {3/2} \\
&\leq C^\prime  |u(t)|^{6} + \frac{1}{4} \left |\nabla u (t)\right | ^ {2},
 \end{split}
\end{equation*}
\begin{equation}\label{uux1}
\begin{split}
 \epsilon   \left |\int _{\dom}   u_{xxx} \,   u_x \,u\,   d\, \dom  \right|& \leq  \epsilon  \left | u_{xxx}\right| \left | u_x \,u \right|\\
& \leq  C^\prime    \epsilon \left | u_x \,u \right|^2   +  \frac{\epsilon}{10} | u_{xxx}|^2  \\
& \leq  C^\prime   \epsilon |  u   |^2_{L^4 (\dom)}|u _x|^2_{L^4 (\dom)} +    \frac{ \epsilon}{10} | u_{xxx}|^2 \\
&\leq (\mbox{by $H^{3/4}(\dom) \subset L^4 (\dom)$ in 3D}) \\
&\leq C^\prime   \epsilon  |u |^{1/2} |\nabla u |^{3/2}  |u_x | ^{1/2} |u_x |_{H^1 (\dom)}^{3/2} + \frac{\epsilon}{10} | u_{xxx}|^2 \\
& \leq C^\prime   \epsilon  |u |^{1/2} |\nabla u |^{2}  |u |_{H^2 (\dom)}^{3/2} + \frac{\epsilon}{10} | u_{xxx}|^2,
\end{split}
\end{equation}
\begin{align*}
 \epsilon \left |\int _{\dom}   u_{yyy} \,x\,   u_y \,u\,   d\, \dom  \right|& \leq (\mbox{by similar estimates as above})\\
 &\leq  C^\prime  \epsilon  |u |^{1/2} |\nabla u|^{2}  |u |_{H^2 (\dom)}^{3/2} + \frac{\epsilon}{10} | u_{yyy}|^2 ,
\end{align*}
\begin{align*}
 \epsilon \left |\int _{\dom}   u_{zzz} \, x\,  u_z \,u\,   d\, \dom  \right|& \leq (\mbox{by similar estimates as above})\\
 & \leq   C^\prime  \epsilon  |u |^{1/2} |\nabla u |^{2}  |u |_{H^2 (\dom)}^{3/2} + \frac{\epsilon}{10} | u_{zzz}|^2,
\end{align*}
\begin{align*}
\left | \int_{\dom} \Delta f \,   u \, d \dom \right |  &\leq   \left  |\Delta f  \right |^2   +     \left  |  u \right |^2    ,
\end{align*}
\begin{align*}
\left |  \int_{\dom}  f   u ^2 \, d \dom \right |&\leq |f|_{L^{\infty}(\dom)} |u|^2 \leq  |f|^2_{L^{\infty(\dom)}} + |u|^4.
\end{align*}
Collecting the above estimates, along with (\ref{key4})  we observe that
the terms with third-order derivatives in the RHS of  (\ref{uux1}) and the following two inequalities can be canceled by a term  on the LHS of (\ref{key4}). Thus (\ref{key4})  now yields
\begin{equation}\label{key5}
\begin{split}
\frac{1}{4}  |  \nabla u^\epsilon (t) |^2  & + \frac{\epsilon}{10}   \int^t _0 [\nabla u^\epsilon]^2_2 \, ds\\
&\leq      \int^t _0  \left( 1+   C^\prime  \epsilon  |u^\epsilon|^{1/2} |u^\epsilon|_{H^2 (\dom)}^{3/2}  \right)    |\nabla u^\epsilon (s) | ^ {2}\, ds  +   C^\prime  |u^\epsilon(t)|^{6}  +  \kappa_0 \\
&\,\,\,\,\,\,\,
 +  \int^t _0    |\Delta f |^2 \,  ds +  \int^t _0 |u^\epsilon  |^2 + |u^\epsilon|^4 \,   ds + \int^t_0 |f|^2 _{L^{\infty}(\dom)}\, ds\\
& \leq (\mbox{thanks to (\ref{eq10-9})}) \\
 &\leq \int^t _0  \left( 1+   C^\prime  \epsilon  \mu_1^{1/4} |u^\epsilon|_{H^2 (\dom)}^{3/2}  \right)    |\nabla u^\epsilon (s) | ^ {2}\, ds +   C^\prime  \mu_1^{3}  +  \kappa_0\\
& \,\,\,\,\,\,\,+  |f|^2 _{L^2 (0, T; H_0^2 (\dom))} +  ( \mu_1 + \mu_1 ^2 )T + |f|^2 _{L^2 (0, T; L^\infty(\dom))}.
\end{split}
\end{equation}
In particular,   setting  $\sigma ^\epsilon(t) :=  1+   C^\prime  \epsilon  \mu_1^{1/4} |u^\epsilon|_{H^2 (\dom)}^{3/2} $,
   from (\ref{key5}) we  deduce
\begin{equation}\label{key55}
\begin{split}
\frac{1}{4}  |  \nabla u^\epsilon (t) |^2 + \frac{\epsilon}{10}   \int^t _0 [\nabla u^\epsilon]^2_2 \, ds
 &\leq    \int^t _0 \sigma^\epsilon (s) |\nabla u^\epsilon (s) | ^ {2}\, ds  \\
&\,\,\,\,\,\,\,\,\quad+   C^\prime  \mu_1^{3}  +  \kappa_0
 +  |f|^2 _{L^2 (0, T; H_0^2(\dom))}\\
 & \,\,\,\,\,\,\,\, \quad+  ( \mu_1 + \mu_1 ^2 )T + |f|^2 _{L^2 (0, T; L^\infty(\dom))}.
\end{split}
\end{equation}
Since $|u^\epsilon|_{H^2 (\dom)}^{3/2} \leq |u^\epsilon|_{H^2 (\dom)}^{2} + C^\prime  $, we find
\begin{equation*}\label{phi1}
\begin{split}
\int^{T}_0  \sigma^\epsilon (s)  \,ds    &\leq  T +   C^\prime  \epsilon \mu_1^{1/4}  \int ^T_0 \left (|u^\epsilon|_{H^2 (\dom)}^{2} +C^\prime \right ) ds \\
& \leq (\mbox{thanks to (\ref{set1p})}) \\
&  \leq \mbox{const} := \mu_3. \\
\end{split}
\end{equation*}
We can then apply  the Gronwall inequality to (\ref{key55}) to obtain
\begin{equation}\label{eq10-99}
\sup_{t\in (0, T)} |  \nabla u ^\epsilon(t) |^2  + \frac{\epsilon}{10}   \int^T _0 [\nabla u^\epsilon]^2_2 \, ds \leq const:=\mu_4.
\end{equation}
This  together with (\ref{set1}) implies  (\ref{main8}) and (\ref{main8p}).  \qed


\subsubsection{Estimates independent of $\epsilon$ for  $u^\epsilon_{xxx}$ and $u^\epsilon u^\epsilon_x $}
\label{other}
For the sake of the passage to the limit on the boundary conditions  and the compactness argument, we now derive   bounds  independent of $\epsilon$  for $  u^\epsilon_{xxx}$ and  $u^\epsilon\,u^\epsilon_x$. In particular, to obtain the estimates for $u^\epsilon_{xxx}$, we first deduce a bound independent of $\epsilon $ for $\epsilon \, u_{xxxx}^\epsilon$ in $L^{2} (0, T ; L^{2} (\dom))  $.
\begin{prop}\label{localbound2}
Under the same assumptions as in Proposition \ref{localbounds1}, we further suppose that
 \begin{equation}\label{u0xx}
u_{0xx}\in L^2 (\dom),
\end{equation}
\begin{equation}\label{fxxx}
f_{xxx}\in L^2 (0, T; L^2   (\dom)),
\end{equation}
and $f_{xx}$ assume  the periodic boundary condition on $x=0$, $1$.
Then we have the following bounds independent of $\epsilon$,
\begin{equation}\label{uxxxxre}
\epsilon \, [u_{xx}]_2  \mbox{ is  bounded   in } L^{2} (0, T ; L^{2} (\dom))  ,
\end{equation}
\begin{equation}\label{h2}
\begin{split}
u^\epsilon u^\epsilon_x \mbox{ is bounded     in } L^{\infty} (0, T ; L^{3/2} (\dom)) .
\end{split}
\end{equation}
\begin{equation}\label{uxxxre}
 u_{xxx}^\epsilon  \mbox{ is  bounded   in } L^{3/2}(I_x;  H^{-1}_t (0, T;\, H^{-4}( I_{x^\perp}   )) )  ,
\end{equation}
\end{prop}

\noindent \textbf{Proof.} For notational simplicity, we will drop the super index $\epsilon$ in the calculations.  Multiplying (\ref{10-11})    by   $   u^\epsilon_{xxxx} $,   integrating over $\dom$ and integrating by parts we find:
\begin{equation*}\displaystyle
\begin{split}
& \bullet \int_{\dom} u_{t} \,  u_{xxxx} \, d\dom   =(\mbox{thanks to (\ref{dirch}) and (\ref{uxxx})} )  =  \frac{1}{2} \frac{d}{dt}  |     u_{xx}   |^2 , \\
&  \bullet \int_{\dom} \Delta  u_{x} \,  u_{xxxx}\, d\dom  =(\mbox{thanks to (\ref{dirch}), (\ref{uxxx}) and  (\ref{eq107})} )=0,\\
&  \bullet \int_{\dom}  u_x   u_{xxxx}\, d\dom  =(\mbox{thanks to (\ref{dirch}),(\ref{uxxx}) and  (\ref{eq107})} ) =0, \\
&  \bullet \int_{\dom} u u_{x} \,   u_{xxxx}\, d\dom = - \int_{\dom}  u_{x}^2 \,   u_{xxx}\, d\dom -  \int_{\dom} u u_{xx} \,   u_{xxx}\, d\dom  \\
& \quad\quad\quad\quad  = \frac{5}{2}\int_{\dom}  u_{x} \,   u^2_{xx}\, d\dom,  \\
& \bullet\,\,      \int _{\dom}   u_{yyyy}   u_{xxxx} \,d\, \dom =  (\mbox{thanks to (\ref{dirch}), (\ref{uxxx}) and  (\ref{eq107})} ) \\
 &\quad\quad\quad\quad =\int _{\dom}   u^2_{xxyy}    \,d\, \dom ,  \\
 & \bullet\,\,      \int _{\dom}   u_{zzzz}   u_{xxxx} \,d\,\dom  = \int _{\dom}     u_{xxzz} ^2  \,d\, \dom ,\\
    & \bullet\,\, \int _{\dom } f    u_{xxxx} \, d\dom = - \int _{\dom } f_{xxx} \,   u_{x } \, d\dom \leq |f_{xxx}|^2 + |u_x|^2 .
  \end{split}
\end{equation*}
Hence we find
\begin{equation*}\label{main1p}
\begin{split}
      \frac{1}{2} \frac{d}{dt}  |     u^\epsilon _{xx}   |^2   +   \epsilon [u^\epsilon _{xx}]^2_2  \leq   \frac{5}{2}\int_{\dom}  u ^\epsilon_{x} \,   \left(u^\epsilon_{xx} \right)^2\, d\dom      +    |f_{xxx}|^2 + |u^\epsilon _x|^2.
\end{split} \end{equation*}
Multiplying both sides by $\epsilon$ we obtain
\begin{equation}\label{main1p}
\begin{split}
      \frac{\epsilon}{2} \frac{d}{dt}  |     u^\epsilon _{xx}   |^2   +   \epsilon^2 [u^\epsilon _{xx}]^2_2  \leq   \frac{5\epsilon}{2}\int_{\dom}  u ^\epsilon _{x} \,  \left(  u^\epsilon _{xx} \right )^2\, d\dom   + \epsilon   |f_{xxx}|^2 +\epsilon |u^\epsilon _x|^2.
\end{split} \end{equation}
We estimate the first term on the right-hand side of (\ref{main1p}) and find
\begin{align*}
\epsilon \left | \int_{\dom}  u_{x} \,   u^2_{xx}\, d\dom \right | &\leq \epsilon\,|u_x| |u_{xx}| ^2 _{L^4 (\dom)}\\
& \leq C^\prime \epsilon\, |u_x||u_{xx}|^{1/2} |\nabla u_{xx}|^{3/2}\\
 & \leq (\mbox{by the intermediate derivative theorem } \left|u_{xx}\right| ^2 \leq \left |u_{x}\right|  \left|u_{xxx}\right|)\\
& \leq C^\prime \epsilon |u_{x}|^{5/4}  |u_{xxx}|^{1/4} |\nabla u_{xx}|^{3/2}\\
& \leq C^\prime  \epsilon   |u_{x}|^{5/4} |\nabla u_{xx}|^{7/4}\\
 & \leq (\mbox{thanks to (\ref{eq10-99})})\\
 & \leq C^\prime  \epsilon \mu_4^{5/8} |\nabla u_{xx}|^{7/4}.
\end{align*}
This  along with (\ref{main1p}) implies
\begin{equation*}
\begin{split}
      \frac{\epsilon}{2} \frac{d}{dt}  |  u^\epsilon _{xx}   |^2   +   \epsilon^2 [u^\epsilon _{xx}]^2_2  \leq  C^\prime \epsilon \, \mu_4^{5/8} |\nabla u^\epsilon_{xx}|^{7/4}+ \epsilon   |f_{xxx}|^2 +\epsilon \mu_4 .
\end{split} \end{equation*}
Integrating both sides in $t$ from $0$ to $T$, we find
 \begin{equation}\label{main11p}
\begin{split}
          \epsilon^2 \int ^T _0 [u^\epsilon _{xx}]^2_2 \,dt \leq  \frac{\epsilon}{2}    |  u _{0xx}   |^2  + C^\prime  \mu_4^{5/8}\int ^T _0\epsilon |\nabla u^\epsilon_{xx}|^{7/4}\, dt+ \epsilon   |f_{xxx}|_{L^2(0, T; L^2 (\dom))}^2 +\epsilon   \mu_4 T.
\end{split} \end{equation}
From (\ref{main8p}), we see that $\int ^T _0\epsilon \left|\nabla u^\epsilon_{xx}\right |^{7/4}\, dt \leq C^\prime \int ^T _0\epsilon \left ( |\nabla u^\epsilon_{xx}|^{2} +1 \right )\, dt  \leq const:=\mu_6$.  This along with (\ref{main11p}) implies (\ref{uxxxxre}).

Now since
\begin{equation*}
\begin{split}
\int_{\dom}( u  u _x)^{3/2} \, d\dom
 \leq C^\prime |u |^{3/2}_{L^6(\dom)} |u _x|^{3/2}
\leq (\mbox{by $H^1(\dom) \subset L^6 (\dom)$ in $3D$})
 \leq C^\prime | u |_{H^1}^3 ,
\end{split}
\end{equation*}
this along with     (\ref{main8}) implies (\ref{h2}), and hence
\begin{equation}\label{h2p}
\begin{split}
u^\epsilon u^\epsilon_x \mbox{ is bounded     in } L^{3/2} (I_x; L^{3/2} ((0, T) \times I_{x^{\perp}})) .
\end{split}
\end{equation}

Finally rewriting (\ref{10-11}) we find
\begin{eqnarray}\label{10-11p}\displaystyle
u^\epsilon _{xxx} = -   u^\epsilon _t- \Delta^\perp   u ^\epsilon _x- c   u^\epsilon _x- u ^\epsilon\, u^\epsilon _x - \epsilon \, u^\epsilon _{xxxx}  -  \epsilon \, u^\epsilon  _{yyyy} - \epsilon \,u ^\epsilon _{zzzz}.
\end{eqnarray}
Thanks to (\ref{uxxxxre}), we see that $\epsilon u^\epsilon_{xxxx}$ remains bounded in $L^2 (0, T; \, L^2 (\dom))$. Moreover since $u^\epsilon $ remains $L^\infty (0, T; H^1 (\dom))$,  we find that each term on the right-hand side of (\ref{10-11p}) except for $u ^\epsilon\, u^\epsilon _x$ remains bounded at least in  $L^{2} (I_x; H^{-1}_t (0, T;\, H^{-4}( I_{x^\perp}   )) )$. This together with (\ref{h2p}) implies that each term on the right-hand side of (\ref{10-11p}) remains bounded at least in  $L^{3/2} (I_x; H^{-1}_t (0, T;\, H^{-4}( I_{x^\perp}   )) )$.
Thus we obtain  (\ref{uxxxre}) from (\ref{10-11p}).

.   \qed

\subsection{The main result}
\label{gob}

Using a compactness argument,
we can pass to the limit in  (\ref{10-11}) and obtain (\ref{eq01}), with a function $u\in  \mathcal C ([0, T];  H^1(\dom))\cap H^3 (I_x;  H^{-1}_t (0, T;\, H^{-4}( I_{x^\perp}   )) )$.
Moreover, from (\ref{uxxxre}) we see that $u^{\epsilon}_{xxx}$ converges weakly in  $L^{3/2}(I_x;  H^{-1}_t (0, T;\, H^{-4}( I_{x^\perp}   )) ) $, hence by the trace theorem and Mazur's theorem, we deduce that $u ^\epsilon_{x^j} (0, x^\perp, t)$ and $u^\epsilon _{x^j} (1, x^\perp, t)$ converge weakly in $ H^{-1}_t (0, T;\, H^{-4}( I_{x^\perp}   ))$, $j=1,\,\,2$. Thus from (\ref{uxxx}) we obtain (\ref{eq4}).

Now we are ready to state the main result of the article  by collecting all the previous   estimates.
\begin{thm}\label{localstrong}
  The assumptions are the same as  in Proposition \ref{localbound2}, that is (\ref{initial2}),  (\ref{f}),   (\ref{f1}),     (\ref{u02})),  (\ref{fxxx}),  (\ref{u0xx}), and $f$ and $f_{x^j}$ assume the periodic boundary conditions on $x=0$, $1$, $j=1, 2$.
  Then
 the initial and boundary value problem for the ZK equation, that is,   (\ref{eq01}), (\ref{dirch})-(\ref{eq107}) and
   (\ref{eq33}), possesses at least  a solution u:
\begin{equation}\label{defs}
u  \in   \mathcal C ([0, T] ; H^1(\dom))\cap  W^{3,\,3/2} (I_x; H^{-1}_t (0, T;\,   H^{-4}( I_{x^\perp}   )) ) .
\end{equation}

\end{thm}
\begin{rem}
\label{uav}
We can   obtain stronger regularity for $\bar u (x^\perp, t) :=  \int ^ 1_0 u (x, x^\perp, t)\, dx$. Integrating (\ref{eq01}) in $x$ from $0$ to $1$, we find  by (\ref{dirch}) and (\ref{eq4})
\begin{eqnarray}\label{baruisf}\displaystyle
&\dfrac{\partial \bar u }{\partial t}   =  \bar f .
\end{eqnarray}
Thus $u= \bar u + v$, where $\bar u$ satisfies (\ref{baruisf}),  and $v$ satisfies $\bar v =0$ and (\ref{defs}).

\end{rem}


\section{Discussions about the uniqueness of solutions.}
\label{unique}

Let $u$ and $v$ be two solutions of  (\ref{eq01})-(\ref{eq107}) and   (\ref{eq33}) and let $w=u-v$. Letting  $\bar w (x^\perp, t) :=  \int ^ 1_0 u (x, x^\perp, t)\, dx$, we see that
$\dfrac{\partial \bar w}{\partial t}   = 0 $ and hence
\begin{equation}
 \bar w (t)=0, \,\,\,\, \forall\,\,\, t\in [0, T] .
\end{equation}
However, it is not clear if we can further prove that $ w (t)=0, \,\,\,\, \forall\,\,\, t\in [0,T]$.
Firstly, the ideas in the proof of existence can not be extended to prove the uniqueness because the structure of the nonlinear term is changed. Secondly,  the methods in \cite{SautTemam} and \cite{SautTemamChuntian} are not applicable  due to the lack of assumptions on the boundary condition $u_x$ at $x=1$.
 For the same reason, the proof of the local existence in \cite{CW} fails as well, which   prevents us from using the methods   in \cite{CaoTiti}.

To conclude,    the uniqueness of solutions in both  dimensions $2$ and $3$ are  still open  due to the partially hyperbolic feature of this model.



\begin{rem}
As for the periodic case, that is, (\ref{eq01}) and the boundary and initial conditions   (\ref{dirch}),  (\ref{eq4}),  (\ref{eq1-26p}) and  (\ref{eq33}), the results are exactly the same as in the Dirichlet case  discussed above. The reasoning is totally the same and therefore we skip it.
\end{rem}

\section*{Acknowledgments}
This work was partially supported by the National Science Foundation under the grants, DMS-0906440  and DMS 1206438, and by the Research Fund of Indiana University.

%
%

The author would like to thank my  advisers Professor Roger Temam and Nathan Glatt-Holtz for their encouragements and suggestions.


   \newpage

\footnotesize
\bibliographystyle{amsalpha}

\bibliography{ref-3}

\normalsize

%

\end{document}